\documentclass[oneside,english]{amsart}
\usepackage[T1]{fontenc}
\usepackage[latin9]{inputenc}
\usepackage{verbatim}
\usepackage{amssymb}

\makeatletter
\numberwithin{equation}{section} 
\numberwithin{figure}{section} 
  \theoremstyle{plain}
  \newtheorem{thm}{Theorem}[section]
  \theoremstyle{remark}
  \newtheorem*{acknowledgement*}{Acknowledgement}
  \theoremstyle{plain}
  \newtheorem{prop}[thm]{Proposition}

\usepackage{alex}
\usepackage[all]{xy}
\renewcommand{\theenumi}{{\it\roman{enumi})}}
\renewcommand{\labelenumi}{\theenumi}
\date{}

\usepackage{babel}
\makeatother

\begin{document}

\title{The smooth representations of $\mathrm{GL}_{2}(\mathcal{O})$}

\author{Alexander Stasinski}

\address{DPMMS, University of Cambridge, Wilberforce Rd, Cambridge, CB3 0WB,
U.~K.}

\email{a.stasinski@dpmms.cam.ac.uk}

\begin{abstract}
We present some unpublished results of Kutzko together with results
of Hill, giving a classification of the smooth (complex) representations
of $\mathrm{GL}_{2}(\mathcal{O})$, where $\mathcal{O}$ is the ring
of integers in a local field with finite residue field.
\end{abstract}
\maketitle

\section{Introduction}

Let $\mathcal{O}$ be the ring of integers in a local field with finite
residue field, and let $\mathfrak{p}$ be its maximal ideal. The smooth
representations of $\mbox{GL}_{2}(\mathcal{O})$, or equivalently,
the representations of the finite groups $\mbox{GL}_{2}(\mathcal{O}/\mathfrak{p}^{r})$,
$r\geq1$, have in one form or another been known to some mathematicians
since the late 70s. There has been two different approaches to this
problem. On the one hand, there is the Weil representation approach
of Nobs and Wolfart (applied to the case $\mathcal{O=\mathbb{Z}}_{p}$
in \cite{MR0499009}). On the other hand, there is the approach via
orbits and Clifford theory due to Kutzko (unpublished), and independently
to Nagornyj \cite{nagornyj1}.

The related case of $\mathrm{SL}_{2}(\mathbb{Z}_{p})$, $p\neq2$,
has been studied by Kloosterman \cite{Kloosterman_I,Kloosterman_II},
Tanaka \cite{Tanaka1,Tanaka2}, Kutzko (unpublished), and Shalika
(for general $\mathcal{O}$ and $p\neq2$) \cite{Shalika}. Another
description of the representations of $\mathrm{SL}_{2}(\mathbb{Z}_{p})$
(including the case $p=2$) was given by Nobs and Wolfart \cite{MR0444787,MR0444788},
using Weil representations. The case $\mathrm{PGL}_{2}(\mathcal{O})$
with odd residue characteristic has been treated by Silberger \cite{Silberger}. 

The series of papers by Hill \cite{Hill_Jordan,Hill_nilpotent,Hill_regular,Hill_semisimple_cuspidal}
use the method of orbits and Clifford theory as a basis for dealing
with the general case of $\mathrm{GL}_{n}(\mathcal{O})$, and contain
several general results, albeit a complete classification only in
certain cases. In particular, Hill's work gives a construction (up
to a description of orbits) of the so called \emph{regular} representations
(resp. \emph{split regular} in the odd conductor case). 

Despite all these works, the literature on the representations of
$\mbox{GL}_{2}(\mathcal{O})$ remains in a somewhat unsatisfactory
state: The Weil representation approach of Nobs and Wolfart works
only for $\mbox{GL}_{2}$ and closely related groups of small rank,
and Nagornyj's paper, in addition to restricting to the case $\mathcal{O}=\mathbb{Z}_{p}$
with $p$ odd, omits certain details and uses an ad hoc argument in
the cuspidal case. Moreover, although Hill's construction of regular
representations covers most of the representations of $\mbox{GL}_{2}(\mathcal{O})$,
it omits a discussion of the orbits, and does not work for cuspidal
representations in the odd conductor case. Hill does have a construction
of so called strongly semisimple representations which include the
cuspidals in the odd conductor case, but this involves a certain implicit
step (see Subsect.~\ref{sub:The-cuspidal-case}), and a more direct
construction of the cuspidal representations is desirable.

The purpose of this paper is to give a complete and uniform construction
of all the irreducible representations of the groups $\mbox{GL}_{2}(\mathcal{O}/\mathfrak{p}^{r})$,
for $r\geq2$ (the case $r=1$ being well-known classically), via
the approach of orbits and Clifford theory. This amounts to a reconstruction
of Kutzko's unpublished results, using whenever possible, the general
results of Hill in order to fit the present construction into a more
general picture.

\begin{acknowledgement*}
The author is grateful to P. Kutzko for lecturing on this material
and answering subsequent questions, and to S. Stevens for many helpful
discussions, in particular for patiently explaining the details of
the cuspidal case. This work was supported by EPSRC Grant GR/T21714/01.
\end{acknowledgement*}

\section{Orbits and Clifford theory}

Let $F$ be a local field with finite residue field $k$ of arbitrary
characteristic. Denote by $\mathcal{O}$ its ring of integers, by
$\mathfrak{p}$ its maximal ideal, and by $\varpi$ a prime element.
For any integer $r\geq1$ we write $\mathcal{O}_{r}$ for the finite
ring $\mathcal{O}/\mathfrak{p}^{r}$. Let $G_{r}=\textrm{GL}_{2}(\mathcal{O}_{r})$,
$K_{i}=\{g\in G_{r}\mid g\equiv1\pmod{\mathfrak{p}^{i}}\}$, for $1\leq i\leq r-1$.
From now on let $i\geq r/2$; then $x\mapsto1+\varpi^{i}x$ induces
an isomorphism\[
M_{2}(\mathcal{O}_{r-i})\longiso K_{i}.\]
The group $G_{r}$ acts on $M_{2}(\mathcal{O}_{r-i})$ by conjugation,
via its quotient $G_{r-i}$. This action is transformed by the above
isomorphism into the action of $G_{r}$ on the normal subgroup $K_{i}$. 

Fix an additive character $\psi:\mathcal{O}\rightarrow\mathbb{C}^{\times}$
with conductor $\mathfrak{p}^{r}$ and define for any $\beta\in M_{2}(\mathcal{O}_{r-i})$,
a character $\psi_{\beta}:K_{i}\rightarrow\mathbb{C}^{\times}$ by\[
\psi_{\beta}(x)=\psi(\Tr(\beta(x-1))).\]
Then $\beta\mapsto\psi_{\beta}$ gives an isomorphism \[
M_{2}(\mathcal{O}_{r-i})\cong\textrm{Hom}(K_{i},\mathbb{C}^{\times}).\]
 For $g\in G_{r}$ we have \[
\psi_{g\beta g^{-1}}(x)=\psi(\Tr(g\beta g^{-1}(x-1)))=\psi(\Tr(\beta g^{-1}(x-1)g))=\psi_{\beta}(g^{-1}xg).\]
Thus, the above isomorphism transforms the action of $G_{r}$ on $M_{2}(\mathcal{O}_{r-i})$
into (the inverse) conjugation of characters.

We will make use of the following well-known results of Clifford theory. 

\begin{thm}
\label{thm:Clifford}Let $G$ be a finite group, and $N$ a normal
subgroup. For any irreducible representation $\rho$ of $N$, define
the \emph{stabilizer}, $T(\rho)=\{g\in G\mid\rho^{g}\cong\rho\}$,
of $\rho$. Then the following hold.
\begin{enumerate}
\item \label{Clifford1}If $\pi$ is an irreducible representation of $G$,
then $\pi|_{N}=e(\bigoplus_{\rho\in\Omega}\rho)$, where $\Omega$
is an orbit of irreducible representations of $N$ under the action
of $G$ by conjugation, and $e$ is a positive integer.
\item \label{Clifford2}Suppose that $\rho$ is an irreducible representation
of $N$, and let\begin{align*}
A & =\{\theta\in\Irr(T(\rho))\mid\langle\Res_{N}^{T(\rho)}(\theta),\rho\rangle\neq0\},\\
B & =\{\pi\in\Irr(G)\mid\langle\Res_{N}^{G}(\pi),\rho\rangle\neq0\}.\end{align*}
Then $\theta\mapsto\Ind_{T(\rho)}^{G}(\theta)$ is a bijection of
$A$ onto $B$.
\item Let $H$ be a subgroup of $G$ containing $N$, and suppose that $\rho$
is an irreducible representation of $N$ which has an extension $\tilde{\rho}$
to $H$ (i.e.~$\tilde{\rho}|_{N}=\rho$). Then \[
\Ind_{N}^{H}(\rho)=\bigoplus_{\chi\in\Irr(H/N)}\chi\tilde{\rho},\]
where each $\chi\tilde{\rho}$ is irreducible, and where we have identified
$\Irr(H/N)$ with representations of $H$ that are trivial on $N$. 
\item \label{Clifford3}If $\rho$ is an irreducible representation of $N$
and $T(\rho)/N$ is cyclic, then there exists an extension of $\rho$
to $T(\rho)$.
\end{enumerate}
\end{thm}
For proofs of the above, see for example \cite{Isaacs}, 6.2, 6.11,
6.17, and 11.22, respectively. The above results \ref{Clifford1}
and \ref{Clifford2} show that in order to obtain a classification
of the representations of $G_{r}$, it is enough to classify the orbits
of irreducible representations $\psi_{\beta}$ of a normal subgroup
$K_{i}$, and to describe all the irreducible representations of the
stabilizers $T(\psi_{\beta})$ which contain $\psi_{\beta}$ (when
restricted to $K_{i}$), that is, to decompose $\Ind_{K_{i}}^{T(\psi_{\beta})}(\psi_{\beta})$
into irreducible representations. This is what we will do in the following.

\subsection{\label{sub:Orbits}Orbits}

Set $l=[\frac{r+1}{2}]$, $l'=[\frac{r}{2}]$; thus $l+l'=r$. Then
$K_{l}$ is the largest abelian group among the kernels $K_{i}$,
so we can describe its irreducible representations $\psi_{\beta}$,
and at the same time it has the advantage of being sufficiently close
to the stabilizers of its representations to enable us to describe
the irreducible components of $\Ind_{K_{i}}^{T(\psi_{\beta})}(\psi_{\beta})$. 

We will consider $G_{r}$-orbits on $M_{2}(\mathcal{O}_{l'})=M_{2}(\mathcal{O}_{r-l})\cong\Hom(K_{l},\mathbb{C}^{\times})$.
The Jordan and rational canonical forms imply that the orbits fall
into four types, according to their reductions mod $\mathfrak{p}$.
The following is a list of the possible orbits in $M_{2}(k)$:\renewcommand{\theenumi}{{(\arabic{enumi})}} 

\renewcommand{\labelenumi}{\theenumi}

\let\oldenumerate=\enumerate \renewcommand{\enumerate}{\oldenumerate\setlength{\itemsep}{0.1 cm}}

\begin{enumerate}
\item \label{Orbit1}$\begin{pmatrix}a & 0\\
0 & a\end{pmatrix}$,
\item $\begin{pmatrix}a & 0\\
0 & d\end{pmatrix}$, where $a\neq d$,
\item \label{Orbit3}$\begin{pmatrix}\textrm{0} & 1\\
-\Delta & s\end{pmatrix}$, where $x^{2}-sx+\Delta$ is irreducible,
\item \label{Orbit4}$\begin{pmatrix}a & 1\\
0 & a\end{pmatrix}$.
\end{enumerate}
If $x\in M_{2}(\mathcal{O}_{i})$, for some $i\geq1$, we shall write
$\overline{x}$ for the image of $x$ in $M_{2}(k)$. Let $\beta\in M_{2}(\mathcal{O}_{l'})$.
Then $\psi_{\beta}|_{K_{r-1}}=\psi_{\overline{\beta}}$, and so if
$\beta'$ is another matrix in $M_{2}(\mathcal{O}_{l'})$, we have
\[
\beta\equiv\beta'\pmod{\mathfrak{p}}\Longleftrightarrow\psi_{\beta}|_{K_{r-1}}=\psi_{\beta'}|_{K_{r-1}}.\]
If $aI_{2}\in M_{2}(k)$ is a scalar matrix ($I_{2}$ is the unit
matrix), then $\psi_{aI_{2}}(x)=\psi(a\Tr(x-1))=\psi(a(\det(x)-1))$.
Thus $\psi_{aI_{2}}=\chi_{a}\circ\det$, for some character $\chi_{a}$
of the group $\{g\in\mathcal{O}_{r}^{\times}\mid g\equiv1\pmod{\mathfrak{p}^{r-1}}\}$.
Since $\mathcal{O}_{r}^{\times}$ is abelian, $\chi_{a}$ has an extension
to $\mathcal{O}_{r}^{\times}$. Consequently $\psi_{aI_{2}}$ has
an extension to $G_{r}$, which we denote by $\tilde{\psi}_{aI_{2}}$.

Suppose now that $\pi$ is a representation of $G_{r}$ such that
$\pi|_{K_{l}}$ contains $\psi_{\beta}$, with $\beta$ a matrix of
type \ref{Orbit1} such that $\overline{\beta}=aI_{2}$. Then $(\pi\otimes\tilde{\psi}_{-aI_{2}})|_{K_{r-1}}=\psi_{0}|_{K_{r-1}}=\mathbf{1}_{K_{r-1}}$.
Thus $\pi$ factors through $G_{r-1}$ after twisting by a one-dimensional
character, i.e.~$\pi$ has the group $K_{r-1}\cap\textrm{SL}_{2}(\mathcal{O}_{r})$
in its kernel. Working inductively, we assume that the representations
of $G_{r-1}$ are known. Hence the representations corresponding to
orbits of type \ref{Orbit1} are described.

Similarly, if $\pi$ is a representation of $G_{r}$ corresponding
to an orbit of type \ref{Orbit4} containing an element $\beta$ such
that $\overline{\beta}=aI_{2}+\left(\begin{smallmatrix}0 & 1\\
0 & 0\end{smallmatrix}\right)$, then $\pi\otimes\tilde{\psi}_{-aI_{2}}$ corresponds to an orbit
in $M_{2}(\mathcal{O}_{l'})$ that can be represented by a matrix
$\beta'$ such that $\beta'\equiv\left(\begin{smallmatrix}0 & 1\\
0 & 0\end{smallmatrix}\right)\pmod{\mathfrak{p}}$. This implies that $\beta'=\left(\begin{smallmatrix}a & 1+b\\
c & d\end{smallmatrix}\right)$ for some $a,b,c,d\in\mathfrak{p}$. Conjugating $\beta'$ by $\left(\begin{smallmatrix}1 & 0\\
a & 1+b\end{smallmatrix}\right)$, we can assume that $\beta'$ has the form $\left(\begin{smallmatrix}0 & 1\\
-\Delta & s\end{smallmatrix}\right)$ for some $\Delta,s\in\mathfrak{p}$. 

If $\beta$ represents an orbit of type \ref{Orbit3}, such that $\beta\equiv\left(\begin{smallmatrix}0 & 1\\
-\Delta & s\end{smallmatrix}\right)\pmod{\mathfrak{p}}$, then $\beta=\left(\begin{smallmatrix}a & 1+b\\
-\Delta+c & s+d\end{smallmatrix}\right)$ for some $a,b,c,d\in\mathfrak{p}$, and conjugating by $\left(\begin{smallmatrix}1 & 0\\
a & 1+b\end{smallmatrix}\right)$, we can assume that $\beta$ has the form $\left(\begin{smallmatrix}0 & 1\\
-\Delta & s\end{smallmatrix}\right)$ for some $\Delta,s\in\mathcal{O}_{l'}.$

We have therefore reduced the problem of describing the representations
of $G_{r}$ to describing the representations corresponding to orbits
in $M_{2}(\mathcal{O}_{l'})$, represented by elements of the following
types\renewcommand{\theenumi}{{(\arabic{enumi}$'$)}} 

\renewcommand{\labelenumi}{\theenumi}

\begin{enumerate}
\item \label{Type1}$\begin{pmatrix}a & 0\\
0 & d\end{pmatrix}$, with $a\not\equiv d\pmod{\mathfrak{p}}$,
\item \label{Type2}$\begin{pmatrix}0 & 1\\
-\Delta & s\end{pmatrix}$, where $x^{2}-sx+\Delta$ is irreducible mod $\mathfrak{p}$,
\item \label{Type3}$\begin{pmatrix}0 & 1\\
-\Delta & s\end{pmatrix}$, where $\Delta,s\in\mathfrak{p}$.
\end{enumerate}
All of these orbits are \emph{regular} in the sense of Hill \cite{Hill_regular}.
Indeed, by Theorem 3.6 in \cite{Hill_regular}, an element $\beta\in M_{2}(\mathcal{O}_{l'})$
is regular if and only if the centralizer $C_{\mathbf{G}_{1}}(\overline{\beta})$
is abelian, where $\mathbf{G}_{1}=\mathrm{GL}_{2}(\overline{k})$,
and $\overline{k}$ is an algebraic closure of $k$. Now simple calculations
show that $C_{\mathbf{G}_{1}}(\overline{\beta})$ is conjugate to
the group \[
\left(\begin{smallmatrix}*\vspace{.1cm}\hspace{.1cm} & 0\\
0\hspace{.1cm} & *\end{smallmatrix}\right),\textrm{ }\left\{ \left(\begin{smallmatrix}a\vspace{.1cm} & b\\
-\overline{\Delta}b\hspace{.1cm} & a+\overline{s}b\end{smallmatrix}\right)\in\textrm{GL}_{2}(\overline{k})\right\} ,\textrm{ or }\left\{ \left(\begin{smallmatrix}a\vspace{.1cm}\hspace{.1cm} & *\\
0\hspace{.1cm} & a\end{smallmatrix}\right)\mid a\in k^{\times}\right\} ,\]
when $\beta$ is of type \ref{Type1}, \ref{Type2}, or \ref{Type3},
respectively.

We observe that the representations corresponding to orbits of type
\ref{Type1} and \ref{Type3} are \emph{split regular} in the sense
of Hill (cf.~\cite{Hill_regular}, 4.4). This means that the orbits
are regular and that the characteristic polynomial splits completely
mod $\mathfrak{p}$ (not necessarily into distinct factors). The representations
corresponding to orbits of type \ref{Type2} are \emph{cuspidal} (cf.~\cite{Hill_semisimple_cuspidal},
4.2), i.e.~the characteristic polynomial is irreducible mod $\mathfrak{p}$.

Let $\beta\in M_{2}(\mathcal{O}_{l'})$ be an element in an orbit
of one of the above three types. Since each of these orbits is regular,
it follows from Corollary 3.7 in \cite{Hill_regular} that the stabilizer
of $\psi_{\beta}$ in $G_{r}$ is given by \[
T(\psi_{\beta})=\mathcal{O}_{r}[\hat{\beta}]^{\times}K_{l'},\]
where $\hat{\beta}$ is a lift of $\beta$ to an element in $M_{2}(\mathcal{O}_{r})$. 

%
{}

To describe the representations of $G_{r}$ corresponding to the above
orbits, we shall proceed by considering two cases depending on whether
$r$ is even or odd.

\section{The case $r=2l,\ l'=l$}

Let $\beta$ be a representative of one of the three orbits above,
and choose a lift $\hat{\beta}\in M_{2}(\mathcal{O}_{r})$. Let $\theta\in\Hom(\mathcal{O}_{r}[\hat{\beta}]^{\times},\mathbb{C}^{\times})$
be a character such that $\theta$ is equal to $\psi_{\beta}$ when
restricted to $K_{l}\cap\mathcal{O}_{r}[\hat{\beta}]^{\times}$ (such
a $\theta$ exists since $\mathcal{O}_{r}[\hat{\beta}]^{\times}$
is abelian). Then $\theta\psi_{\beta}$ is a one-dimensional representation
of $T(\psi_{\beta})$ such that $\theta\psi_{\beta}(xy)=\theta(x)\psi_{\beta}(y)$,
for $x\in\mathcal{O}_{r}[\hat{\beta}]^{\times}$, $y\in K_{l}$.

Define the representation \[
\pi(\theta,\beta)=\Ind_{\mathcal{O}_{r}[\hat{\beta}]^{\times}K_{l}}^{G_{r}}(\theta\psi_{\beta}).\]
By Theorem \ref{thm:Clifford} \ref{Clifford2}, $\pi(\theta,\beta)$
is an irreducible representation, $\pi(\theta,\beta)=\pi(\theta',\beta)\Rightarrow\theta=\theta'$,
and every irreducible representation of $G_{r}$ corresponding to
an orbit of one of the above three types is obtained in this way.
This is a special case of Theorem 4.1 in \cite{Hill_regular}, which
holds for any $\mathrm{GL}_{n}$ and any regular orbit.

We observe that a different choice of representative $\beta$ for
the orbit, gives a conjugate character $\psi_{\beta}$ and a conjugate
stabilizer $T(\psi_{\beta})$. The resulting induced representations
are therefore isomorphic for any choice of representative $\beta$.
Moreover, a different choice of $\hat{\beta}$, even though it may
give a different group $\mathcal{O}_{r}[\hat{\beta}]^{\times}$ and
different characters $\theta$, will give a resulting set of representations
$\theta\psi_{\beta}$, that equals the original one. This is because
in both cases we obtain the set of irreducible components of $\Ind_{K_{l}}^{T(\psi_{\beta})}(\psi_{\beta})$.
Thus, the set of representations $\pi(\theta,\beta)$ depends only
on the orbit of $\beta$, and is independent of the choice of lift
$\hat{\beta}$.

\section{The case $r=2l-1,\ l'=l-1$}

As before, Let $\beta$ be a representative of one of the three orbits
above, and choose a lift $\hat{\beta}\in M_{2}(\mathcal{O}_{r})$.
The thing that makes this case more difficult is that we start from
a character $\psi_{\beta}$ on $K_{l}$, while $T(\psi_{\beta})=\mathcal{O}_{r}[\beta]^{\times}K_{l'}$,
and $K_{l}\lneq K_{l'}$, so it is not as easy to decompose $\Ind_{K_{l}}^{T(\psi_{\beta})}(\psi_{\beta})$.
This the {}``odd conductor case'' in turn falls into two subcases:
the split case, and the cuspidal case.

\subsection{The split case}

Following Hill \cite{Hill_regular}, Sect.~4, we assume that $\beta$
lies in a split regular orbit, i.e.~of type \ref{Type1} or \ref{Type3}.
Let $B_{r}$ be the subgroup of $G_{r}$ of upper-triangular matrices.
By the proof of Lemma 4.5 in \cite{Hill_regular}, the group $H_{\beta}=K_{l}(B_{r}\cap K_{l'})$
has the property that $H_{\beta}/K_{l}$ is a maximal isotropic subspace
of $K_{l'}/K_{l}$ with respect to the form%
\footnote{The form takes values in $k$, not in $\mathbb{C}$ as stated by Hill.%
} \[
\langle\cdot,\cdot\rangle_{\overline{\beta}}:K_{l'}/K_{l}\times K_{l'}/K_{l}\longrightarrow k,\qquad\langle xK_{l},yK_{l}\rangle_{\overline{\beta}}=\Tr(\overline{\beta}(\overline{m}\overline{n}-\overline{n}\overline{m})),\]
where $x=1+\varpi^{l'}m$ and $y=1+\varpi^{l'}n$. Moreover, $H_{\beta}$
is normal in $T(\psi_{\beta})$, and we have the following diagram
of groups. $$\xymatrix{& T(\psi_\beta)\ar@{-}[d] \\ K_{l'}\ar@{-}[ru]\ar@{-}[d] & \mathcal{O}_r[\hat{\beta}]^{\times}H_{\beta}\ar@{-}[d] \\ H_{\beta}\ar@{-}[ru]\ar@{-}[d] & \mathcal{O}_r[\hat{\beta}]^{\times}K_l\\K_l\ar@{-}[ru] &\\}$$We
now review the main steps of the proof of Theorem 4.6 in \cite{Hill_regular}.
Let $N=\ker(\psi_{\beta})$ and write $\overline{\psi}_{\beta}$ for
$\psi_{\beta}$ viewed as a character of $K_{l}/N$. One can show
that $H_{\beta}/N$ and $\mathcal{O}_{r}[\hat{\beta}]^{\times}K_{l}/N$
are both abelian. Hence, $\overline{\psi}_{\beta}$ can be extended
to both these groups. Next it is shown that there exists exactly $q^{2}$
($q^{n}$ in the case of $\textrm{GL}_{n}$) $\mathcal{O}_{r}[\hat{\beta}]^{\times}H_{\beta}/N$-stable
extensions of $\overline{\psi}_{\beta}$ to $H_{\beta}/N$. Let $\psi_{\beta}'$
be one of these. Furthermore, one can construct an extension $\psi_{\beta}''$
of $\overline{\psi}_{\beta}$ to $\mathcal{O}_{r}[\hat{\beta}]^{\times}K_{l}/N$,
such that $\psi_{\beta}''$ and $\psi_{\beta}'$ agree on $(H_{\beta}/N)\cap(\mathcal{O}_{r}[\hat{\beta}]^{\times}K_{l}/N)$.
Then one can check that the map\[
\psi_{\beta}''':\mathcal{O}_{r}[\hat{\beta}]^{\times}H_{\beta}/N\longrightarrow\mathbb{C},\qquad\psi_{\beta}'''(xhN)=\psi_{\beta}''(xN)\psi_{\beta}'(hN),\]
where $x\in\mathcal{O}_{r}[\hat{\beta}]^{\times},h\in H_{\beta}$,
is a well-defined homomorphism. Denote by $\tilde{\psi}_{\beta}$
the composition of $\psi_{\beta}'''$ with the natural homomorphism
$\mathcal{O}_{r}[\hat{\beta}]^{\times}H_{\beta}\rightarrow\mathcal{O}_{r}[\hat{\beta}]^{\times}H_{\beta}/N$.
Thus $\tilde{\psi}_{\beta}$ is an extension of $\psi_{\beta}$ to
$\mathcal{O}_{r}[\hat{\beta}]^{\times}H_{\beta}$.

Now define the characters\[
\rho_{\beta}=\Ind_{H_{\beta}}^{K_{l'}}(\tilde{\psi}_{\beta}|_{H_{\beta}}),\textrm{ and }\zeta_{\beta}=\Ind_{\mathcal{O}_{r}[\hat{\beta}]^{\times}H_{\beta}}^{T(\psi_{\beta})}(\tilde{\psi}_{\beta}).\]
It follows from \cite{Hill_regular}, Prop.~4.2, that $\rho_{\beta}$
is irreducible. Moreover, it is shown that $\zeta_{\beta}$ is an
extension of $\rho_{\beta}$. Frobenius reciprocity and counting degrees
shows that the induced representation $\Ind_{K_{l'}}^{T(\psi_{\beta})}(\rho_{\beta})$
decomposes as the sum of the representations $\omega\zeta_{\beta}$,
where $\omega$ runs through the linear characters of the abelian
group $T(\psi_{\beta})/K_{l'}$. Finally, it is shown that every irreducible
component of $\Ind_{K_{l}}^{K_{l'}}(\psi_{\beta})$ has the form $\Ind_{H_{\beta}}^{K_{l'}}(\psi_{\beta}')$
for some $\mathcal{O}_{r}[\hat{\beta}]^{\times}$-stable extension
$\psi_{\beta}'$ of $\psi_{\beta}$. Thus, every irreducible component
of $\Ind_{K_{l}}^{T(\psi_{\beta})}(\psi_{\beta})$ has the form $\omega\zeta_{\beta}$,
where $\zeta_{\beta}$ is induced from some $\tilde{\psi}_{\beta}$,
and $\omega$ is as above, and all the characters $\omega\zeta_{\beta}$
are distinct for distinct $\psi_{\beta}'$ and $\omega$.

Thus, every representation $\omega\zeta_{\beta}$ is isomorphic to
an induced representation \[
\Ind_{\mathcal{O}_{r}[\hat{\beta}]^{\times}H_{\beta}}^{T(\psi_{\beta})}(\omega'\tilde{\psi}_{\beta}),\]
where $\omega'$ is a linear character of $\mathcal{O}_{r}[\hat{\beta}]^{\times}H_{\beta}/H_{\beta}$,
and distinct $\omega'$ give distinct $\omega\zeta_{\beta}$.

\subsection{\label{sub:The-cuspidal-case}The cuspidal case}

We now assume that $\beta$ is in a cuspidal orbit, i.e.~of type
\ref{Type2}. In \cite{Hill_semisimple_cuspidal}, Prop.~3.6 Hill
gives a construction of so called strongly semisimple representations,
which include the cuspidal ones as a special case. This construction
is less explicit than ours because it involves the non-constructive
existence of an irreducible constituent of $\Ind_{K_{l}}^{K_{l'}}(\psi_{\beta})$.
For the construction of the cuspidal representations, there are also
several slightly different methods giving parametrizations in terms
of certain characters of $\mathcal{O}_{F_{2}}^{\times}$, where $F_{2}$
is the quadratic unramified extension of $F$ (cf.~\cite{AOPS} and
its references). The method we give here is similar to a more general
construction used by Bushnell and Kutzko (cf.~\cite{BushnellKutzko})
in their description of supercuspidal representations of $\mbox{GL}_{n}(F)$.
In particular, the method can be adapted to the construction of cuspidal
representations for any $\mbox{GL}_{n}(\mathcal{O})$. We will make
use of the following result from \cite{Bushnell_Frohlich}, 8.3.3.

\begin{prop}
\label{pro:Bushnell-Frohlich}Let $G$ be a finite group and $N$
a normal subgroup, such that $G/N$ is an elementary abelian $p$-group.
Thus $G/N$ has a structure of $\mathbb{F}_{p}$-vector space. Let
$\chi$ be a one-dimensional representation of $N$, which is stabilized
by $G$. Define a alternating bilinear form\[
h_{\chi}:G/N\times G/N\longrightarrow\mathbb{C}^{\times},\qquad h_{\chi}(g_{1}N,g_{2}N)=\chi([g_{1},g_{2}])=\chi(g_{1}g_{2}g_{1}^{-1}g_{2}^{-1}).\]
Assume that the form $h_{\chi}$ is non-degenerate. Then there exists
a unique up to isomorphism irreducible representation $\rho_{\chi}$
of $G$ such that $\rho_{\chi}|_{N}$ contains $\chi$.
\end{prop}
Note that the form is well-defined%
\footnote{ $\chi([g_{1}n_{1},g_{1}n_{2}])= \chi(g_{1}g_{2}(g_{2}^{-1}n_{1}g_{2}n_{2}n_{1}^{-1}g_{1}^{-1}n_{2}^{-1}g_{1})g_{1}^{-1}g_{2}^{-1})\\ \indent =\chi(g_{1}g_{2}(g_{2}^{-1}n_{1}g_{2}n_{2}n_{1}^{-1}g_{1}^{-1}n_{2}^{-1}g_{1})(g_{1}g_{2})^{-1}[g_{1},g_{2}])\\ \indent =\chi(g_{2}^{-1}n_{1}g_{2}n_{2}n_{1}^{-1}g_{1}^{-1}n_{2}^{-1}g_{1})\chi([g_{1},g_{2}])=\chi(g_{2}^{-1}n_{1}g_{2})\chi(n_{2}n_{1}^{-1})\chi(g_{1}^{-1}n_{2}^{-1}g_{1})\chi([g_{1},g_{2}])\\ \indent =\chi([g_{1},g_{2}]).$%
}, so the non-trivial thing here is the uniqueness part. 

Set $Z^{1}=\mathcal{O}_{r}[\hat{\beta}]^{\times}\cap K_{1}$. The
strategy to find the irreducible constituents of $\Ind_{K_{l}}^{T(\psi_{\beta})}(\psi_{\beta})$
is as follows. First we extend $\psi_{\beta}$ to a representation
$\tilde{\psi}_{\beta}$ of $Z^{1}K_{l}$, as in the even conductor
case. Next we apply the above proposition to the case where $G=Z^{1}K_{l'}$,
$N=Z^{1}K_{l}$, and $\chi=\tilde{\psi}_{\beta}$. We thus obtain
a representation $\eta_{\beta}$ containing $\psi_{\beta}$, and it
is clear that every representation of $Z^{1}K_{l'}$ containing $\psi_{\beta}$
is of this form. Finally, we use Theorem \ref{thm:Clifford}, \ref{Clifford3}
to lift $\eta_{\beta}$ to a representation of $T(\psi_{\beta})$.
Again it is clear that we have in this way obtained every irreducible
component of $\Ind_{K_{l}}^{T(\psi_{\beta})}(\psi_{\beta})$.

To carry out the above strategy, we need to establish several facts
concerning the groups $T(\psi_{\beta})$, $Z^{1}K_{l'}$, $Z^{1}K_{l}$,
and their respective representations: First, to apply Proposition
\ref{pro:Bushnell-Frohlich}, we need show that $Z^{1}K_{l}$ is normal
in $Z^{1}K_{l'}$, that the quotient is elementary abelian, that $Z^{1}K_{l'}$
stabilizes the representation $\tilde{\psi}_{\beta}$, and that the
form $h_{\tilde{\psi}_{\beta}}$ is non-degenerate. Then, to apply
Theorem \ref{thm:Clifford}, \ref{Clifford3}, we need to show that
$Z^{1}K_{l'}$ is normal in $T(\psi_{\beta})$, that $T(\psi_{\beta})$
stabilizes the representation $\eta_{\beta}$, and that the quotient
$T(\psi_{\beta})/Z^{1}K_{l'}$ is cyclic.

We prove all the above facts in one go. First we show that both $Z^{1}K_{l'}$
and $Z^{1}K_{l}$ are normal in $T(\psi_{\beta})$ (and thus $Z^{1}K_{l}$
is normal in $Z^{1}K_{l'}$). It is clear that $Z^{1}K_{l'}$ is normal
in $T(\psi_{\beta})$ since $Z^{1}K_{l'}=T(\psi_{\beta})\cap K_{1}$,
and $K_{1}$ is normal in $G_{r}$. Now let $g\in T(\psi_{\beta})$
and $x\in Z^{1}K_{l}$. Write $g=ek$ with $e\in\mathcal{O}_{r}[\hat{\beta}]^{\times}$
and $k\in K_{l'}$, and write $x=zh$, with $z\in Z^{1}$ and $h\in K_{l}$.
Then \[
gxg^{-1}=z\cdot e([z^{-1},k](khk^{-1}))e^{-1},\]
and the two factors in the bracket lie in $K_{l}$ since $[K_{1},K_{l'}]\leq K_{l}$.
Since $K_{l}$ is normal in $G_{r}$, the whole expression $e([z^{-1},k](khk^{-1}))e^{-1}$
in fact lies in $K_{l}$, so $gxg^{-1}\in Z^{1}K_{l}$, as required.
Now\begin{align*}
\tilde{\psi}_{\beta}(gxg^{-1}) & =\tilde{\psi}_{\beta}(z)\psi_{\beta}(e([z^{-1},k](khk^{-1}))e^{-1})\\
 & =\tilde{\psi}_{\beta}(z)\psi_{\beta}([z^{-1},k])\psi_{\beta}(h)\\
 & =\tilde{\psi}_{\beta}(x)\psi_{\beta}([z^{-1},k]),\end{align*}
where the second equality follows since $T(\psi_{\beta})$ stabilizes
$\psi_{\beta}$. Thus, to show that $T(\psi_{\beta})$ stabilizes
$\tilde{\psi}_{\beta}$, we need only to show that $\psi_{\beta}([z^{-1},k])=1$.
For this, write $k=1+y$, with $y\in\mathfrak{p}^{l'}M_{2}(\mathcal{O}_{r})$
and notice that\[
[z^{-1},k]=z(1+y)z^{-1}(1-y+y^{2})=z(1+y)z^{-1}-y,\]
(here we have used the fact that $y^{3}=0$ and $zy^{2}=y^{2}$ in
$M_{2}(\mathcal{O}_{r})$). Then \[
\psi_{\beta}([z^{-1},k])=\psi(\Tr(\beta(zyz^{-1}-y)))=\psi(\Tr(z(\beta y)z^{-1}-\beta y))=1,\]
where we have used $z\hat{\beta}=\hat{\beta}z\Rightarrow z\beta=\beta z$,
and the fact that $\Tr$ is a conjugacy class function. 

Next, note that since $[K_{l'},K_{l'}]\leq K_{l}$, the bilinear form
does not depend on the lift $\tilde{\psi}_{\beta}$. The form $h_{\psi_{\beta}}=h_{\tilde{\psi}_{\beta}}$
is related to the above form $\langle\cdot,\cdot\rangle_{\overline{\beta}}$
in the following way. First, there is a natural isomorphism\[
\frac{K_{l'}}{(Z^{1}\cap K_{l'})K_{l}}\longiso\frac{Z^{1}K_{l'}}{Z^{1}K_{l}},\qquad k(Z^{1}\cap K_{l'})K_{l}\longmapsto kZ^{1}K_{l},\textrm{ where }k\in K_{l'},\]
and moreover, it is shown in the proof of Prop.~4.2 of \cite{Hill_regular}
that $\langle xK_{l},yK_{l}\rangle_{\overline{\beta}}=0$ if and only
if $\psi_{\beta}([x,y])=1$, so $h_{\psi_{\beta}}(xZ^{1}K_{l},yZ^{1}K_{l})=0$
if and only if $\langle xK_{l},yK_{l}\rangle_{\overline{\beta}}=0$.
By \cite{Hill_regular}, Corollary 4.3%
\footnote{Note that the second equality in the corollary is incorrect.%
} the radical of the form $\langle\cdot,\cdot\rangle_{\overline{\beta}}$
is equal to $(\mathcal{O}_{r}[\hat{\beta}]^{\times}\cap K_{l'})K_{l}/K_{l}$,
so by the above isomorphism the radical of $h_{\psi_{\beta}}$ is
$Z^{1}K_{l}(\mathcal{O}_{r}[\hat{\beta}]^{\times}\cap K_{l'})/Z^{1}K_{l}=Z^{1}K_{l}/Z^{1}K_{l}=1$.
We thus conclude that the form $h_{\psi_{\beta}}$ is non-degenerate.
Since $K_{l'}/K_{l}\cong M_{2}(k)$ is elementary abelian, the same
is true for any of its quotients.

Finally, we show that $T(\psi_{\beta})/Z^{1}K_{l'}$ is cyclic, and
that $T(\psi_{\beta})$ stabilizes the representation $\eta_{\beta}$.
The former is true because the quotient is isomorphic to $\mathcal{O}_{r}[\hat{\beta}]^{\times}/Z^{1}\cong k[\overline{\beta}]^{\times}$,
and this is the multiplicative group of a finite field (the residue
field of the unramified extension generated by a lift of $\beta$
to $M_{2}(\mathcal{O}_{r})$). The latter is true because for $g\in T(\psi_{\beta})$,
the conjugate representation $\eta_{\beta}^{g}$ is another irreducible
representation of $Z^{1}K_{l'}$ whose restriction to $Z^{1}K_{l}$
contains $\tilde{\psi}_{\beta}^{g}=\tilde{\psi}_{\beta}$, so the
uniqueness in Prop.~\ref{pro:Bushnell-Frohlich} implies that $\eta_{\beta}^{g}$
is equivalent to $\eta_{\beta}$.

To conclude, we have shown that any irreducible component of $\Ind_{K_{l}}^{T(\psi_{\beta})}(\psi_{\beta})$
is of the form $\tilde{\eta}_{\beta}$, for some lift $\tilde{\eta}_{\beta}$
of $\eta_{\beta}$ to $T(\psi_{\beta})$. It follows that any irreducible
representation of $G_{r}$ containing $\psi_{\beta}$ is of the form
\[
\Ind_{T(\psi_{\beta})}^{G_{r}}(\tilde{\eta}_{\beta}),\]
for some $\tilde{\eta}_{\beta}$, and distinct inducing representations
give distinct induced representations. Again the dependence on the
choice of lift $\hat{\beta}$ is irrelevant, since for any choice
we obtain the set of irreducible components of $\Ind_{K_{l}}^{T(\psi_{\beta})}(\psi_{\beta})$.
The choice of representative $\beta$ is irrelevant for the same reason
as before.

\bibliographystyle{alex}
\bibliography{alex,Reps_over_fin_rings}

\end{document}